\documentclass[11pt]{amsart}
\usepackage{graphicx}
\usepackage{amssymb, amsmath}
\newcommand{\HG}{H_d(G, \chi)}
\newcommand{\C}{\mathbb{C}}
\begin{document}
\title{Relative symmetric polynomials and money change problem}
\author{\sc M. Shahryari}
\thanks{}
\address{ Department of Pure Mathematics,  Faculty of Mathematical
Sciences, University of Tabriz, Tabriz, Iran }

\email{mshahryari@tabrizu.ac.ir}
\date{\today}

%%% ----------------------------------------------------------------------
\begin{abstract}
This article is devoted to the number of non-negative solutions of the linear Diophantine equation
$$
a_1t_1+a_2t_2+\cdots a_nt_n=d,
$$
where $a_1, \ldots, a_n$, and $d$ are positive integers. We  obtain
a relation between the number of solutions of this equation and
characters of the symmetric group, using {\em relative symmetric
polynomials}. As an application, we  give a necessary and sufficient
condition for the space of the relative symmetric polynomials to be
non-zero.
\end{abstract}

\maketitle


%1111111111111111111111111111111111111111111111111111111111111111111111111
%%%%%%%%%%%%%%%%%%%%%%%%%%%%%%%%%%%%%%%%%%%%%%%%%%%%%%%%%%%%%%%%%%%%%%

\vspace{2cm}

Suppose $a_1, \ldots, a_n$, and $d$ are positive integers, and consider the following linear Diophantine equation:
$$
a_1t_1+a_2t_2+\cdots a_nt_n=d.
$$
Let $Q_d(a_1, \ldots, a_n)$ be the number of non-negative integer solutions of this equation. Computing the exact values of the function $Q_d$
 is the  well-known {\em money change problem}. It is easy to see that a generating function for $Q_d(a_1, \ldots, a_n)$ is
$$
\prod_{i=1}^{n}\frac{1}{1-t^{a_i}}.
$$
This article is devoted for an interesting relation between $Q_d$
and irreducible complex characters of the symmetric group $S_m$,
where $m=a_1+\cdots+a_n$. In fact, we will show that $Q_d$ is a
permutation character of $S_m$, and then we will find its
irreducible constituents. Our main tool, in the investigation of
$Q_d$, is the notion of {\em relative symmetric polynomials}, which
is introduced by the author in \cite{relative}. Once, we find the
irreducible constituents of $Q_d$, we can also obtain a necessary
and sufficient condition for vanishing of the space of relative
symmetric polynomials. A similar result was obtained  in
\cite{permutation}, for vanishing of {\em symmetry classes of
tensors}, using the same method.

We need a survey of results about relative symmetric polynomials in this article. For a detailed exposition, one can see \cite{relative}.

Let $G$ be a subgroup of the full symmetric group $S_m$ of degree $m$ and suppose
$\chi$ is an irreducible complex character of $G$. Let $H_d[x_1,\ldots, x_m]$ be the complex space of
homogenous polynomials of degree $d$ with the independent commuting variables
$x_1, \ldots, x_m$. Suppose $\Gamma^+_{m,d}$ is the set of all $m$-tuples of non-negative
integers, $\alpha=(\alpha_1, \ldots , \alpha_m)$, such that $\sum_{i}\alpha_i=d$. For any $\alpha \in \Gamma^+_{m,d}$, define  $X^{\alpha}$
to be the monomial
$x_1^{\alpha_1}\cdots x_m^{\alpha_m}$. So the set $\{ X^{\alpha}: \alpha \in \Gamma^+_{m,d}\}$ is a basis of $H_d[x_1,\ldots, x_m]$.
We define also an inner product on $H_d[x_1,\ldots, x_m]$ by
$$
<X^{\alpha}, X^{\beta}>=\delta_{\alpha, \beta}.
$$
The group $G$ acts on $H_d[x_1,\ldots, x_m]$ via
$$
q^{\sigma}(x_1, \ldots, x_m)=q(x_{\sigma^{-1}(1)},\ldots, x_{\sigma^{-1}(m)}).
$$
It also acts on $\Gamma^+_{m,d}$ by
$$
\sigma\alpha=(\alpha_{\sigma(1)}, \ldots ,
\alpha_{\sigma(m)}).
$$
Let $\Delta$ be a set of representatives of orbits of $\Gamma^+_{m,d}$ under the action of $G$.

Now consider the idempotent
$$
T(G,\chi)=\frac{\chi(1)}{|G|}\sum_{\sigma\in
G}\chi(\sigma)\sigma
$$
in the group algebra $\C G$. Define the space of {\em relative symmetric polynomials of degree $d$}
with respect to $G$ and $\chi$ to be
$$
\HG=T(G, \chi)(H_d[x_1,\ldots, x_m]).
$$
Let $q\in H_d[x_1,\ldots, x_m]$. Then we set
$$
q^{\ast}=T(G, \chi)(q)
$$
and we call it a {\em symmetrized polynomial} with respect to $G$ and $\chi$.  Clearly
$$
\HG=< X^{\alpha, \ast}: \alpha \in \Gamma^+_{m,d}>,
$$
where $< set\ \  of\ \  vectors >$ denotes the subspace generated by a given set of vectors.

Recall that the inner product of two characters of an arbitrary group $K$ is defined as follows,
$$
[\phi, \psi]_K=\frac{1}{|K|}\sum_{\sigma\in K}\phi(\sigma)\psi(\sigma^{-1}).
$$
In the special case where $K$ is a subgroup of $G$ and $\phi$ and $\psi$ are characters of $G$, the notation $[\phi, \psi]_K$ will
denote the inner product of the restrictions of $\phi$ and $\psi$ to $K$.

It is proved it \cite{relative} that for any $\alpha$, we have
$$
||X^{\alpha, \ast}||^2=\chi(1)\frac{[\chi, 1]_{G_{\alpha}}}{[G:G_{\alpha}]},
$$
where $G_{\alpha}$ is the stabilizer subgroup of $\alpha$ under the action of $G$.
Hence, $X^{\alpha, \ast}\neq 0$, if and only if $[\chi, 1]_{G_{\alpha}}\neq 0$.
According to this result, let $\Omega$ be the set of all $\alpha\in \Gamma^+_{m,d}$, with $[\chi, 1]_{G_{\alpha}}\neq 0$ and suppose
$\bar{\Delta}=\Delta\cap \Omega$.

We proved in \cite{relative}, the following formula for the dimension of $\HG$
$$
\dim \HG= \chi(1)\sum_{\alpha\in \bar{\Delta}}[\chi, 1]_{G_{\alpha}}.
$$
Note that, $\bar{\Delta}$ depends on $\chi$, but $\Delta$ depends only on $G$. Since, $[\chi, 1]_{G_{\alpha}}= 0$,
for all $\alpha\in \Delta-\bar{\Delta}$, we can re-write the above formula, as
$$
\dim \HG= \chi(1)\sum_{\alpha\in \Delta}[\chi, 1]_{G_{\alpha}}.
$$

There is also another interesting formula for the dimension of $\HG$. This is the formula which employs the function $Q_d$ and so it
 connects the money change problem to relative symmetric polynomials. Let
$\sigma \in G$ be any element with the cycle structure
$[a_1,\ldots,a_n]$, ( i.e. $\sigma$ is equal to a product of $n$
disjoint cycles of length $a_1, \ldots, a_n$, respectively). Define
$Q_d(\sigma)$ to be the number of non-negative integer solutions of
the equation
$$
a_1t_1+a_2t_2+\cdots+a_nt_n=d,
$$
so, we have $Q_d(\sigma)=Q_d(a_1, \ldots, a_n)$. If we consider the
free vector space $\mathbb{C}[\Gamma^+_{m,d}]$ as a
$\mathbb{C}G$-module, then for all $\sigma\in G$, we have
$$
Tr\ \sigma=Q_d(\sigma),
$$
and hence, $Q_d$ is a permutation character of $G$.
It is proved in \cite{relative}, that we have also
$$
\dim \HG=\frac{\chi(1)}{|G|}\sum_{\sigma \in G}
\chi(\sigma)Q_d(\sigma).
$$
Note that, we can write this result as $\dim\HG=\chi(1)[\chi, Q_d]_G$. Now, comparing two formulae for the dimension of
$\HG$ and using the {\em reciprocity relation} for induced characters, we obtain
$$
Q_d(a_1, \ldots, a_n)=\sum_{\alpha\in \Delta}(1_{G_{\alpha}})^G(\sigma),
$$
where $\sigma \in S_m$ is any permutation of the cycle structure $[a_1, \ldots, a_n]$, $G$ is any subgroup of $S_m$ containing $\sigma$ and
$m=a_1+\cdots+a_n$. It is clear that, if $\alpha$ and $\beta$ are in the same orbit of $\Gamma^+_{m,d}$, then $(1_{G_{\alpha}})^G=(1_{G_{\beta}})^G$,
so we have also
$$
Q_d(a_1, \ldots, a_n)=\frac{1}{|G|}\sum_{\alpha\in \Gamma^+_{m,d}}|G_{\alpha}|(1_{G_{\alpha}})^G(\sigma).
$$
As our main result in this section, we have,\\

{\bf Theorem A}.
$$
Q_d=\frac{1}{|G|}\sum_{\alpha\in \Gamma^+_{m,d}}|G_{\alpha}|(1_{G_{\alpha}})^G.
$$
\\

In the remaining part of this article, we will assume that, $G=S_m$, then using representation theory of symmetric groups,
we will find, irreducible constituents of $Q_d$.

We need some standard notions from representation theory of symmetric groups.
Ordinary representations of $S_m$ are in
one to one correspondence with {\em partitions} of $m$. Let
$$
\pi=(\pi_1, \ldots, \pi_s)
$$
be any partition of $m$. The irreducible character of $S_m$, corresponding to a partition $\pi$ is denoted by $\chi^{\pi}$. There is also
a subgroup of $S_m$, associated to $\pi$, which is called the {\em Young subgroup} and it is defined as,
$$
S_{\pi}=S_{\{ 1,\ldots, \pi_1\}}\times S_{\{ \pi_1+1,
\ldots, \pi_1+\pi_2\}}\times \cdots.
$$
Therefore, we have $S_{\pi}\cong S_{\pi_1}\times\cdots \times
S_{\pi_s}$.

Let $\pi=(\pi_1, \ldots, \pi_s)$ and $\mu=(\mu_1, \ldots,
\mu_l)$ be two partitions of $m$. We say that $\mu$ majorizes
$\pi$, iff for any $ 1\leq i \leq min\{ s, l\}$, the inequality
$$
\pi_1+\cdots +\pi_i\leq \mu_1+\cdots +\mu_i
$$
holds. In this case we write $\lambda\unlhd\mu$. This is clearly a
partial ordering on the set of all partitions of $m$. A generalized
$\mu$-tableau of type $\pi$ is a function
$$
T:\{ (i,j): 1\leq i\leq h(\mu), 1\leq j\leq \mu_i\} \rightarrow\{ 1,
2, \ldots, m\}
$$
such that for any $1\leq i\leq m$, we have $|t^{-1}(i)|=\pi_i$.
This generalized tableau is called semi-standard if for each $i$,
$j_1<j_2$ implies $T(i,j_1)\leq T(i,j_2)$ and for any $j$, $i_1<i_2$
implies $T(i_1,j)<T(i_2,j)$. In other words, $T$ is semi-standard,
iff every row of $T$ is non-descending and every column of $T$ is
ascending. The number of all such semi-standard tableaux is denoted
by $K_{\mu\pi}$ and it is called the {\em Kostka number}. It is well
known that $K_{\mu\pi}\neq 0$ iff $\mu$ majorizes $\pi$, see
\cite{sagan}, for example. We have also,
\begin{eqnarray*}
K_{\mu\pi}&=&[(1_{S_{\pi}})^{S_m}, \chi^{\mu}]_{S_m}\\
          &=&[1, \chi^{\mu}]_{S_{\pi}}.
\end{eqnarray*}
For any $\alpha\in \Gamma^+_{m,d}$, the multiplicity partition is denoted by $M(\alpha)$, so, to obtain $M(\alpha)$,
we must arrange the multiplicities
of numbers $0, 1, \ldots, d$ in $\alpha$ in the descending order. It is clear that $(S_m)_{\alpha}\cong S_{M(\alpha)}$,
the Young subgroup. If $M(\alpha)=(k_1, \ldots, k_s)$, then we have
$$
|(S_m)_{\alpha}|=k_1!k_2!\ldots k_s!.
$$
In what follows, we use the notation $M(\alpha)!$ for the number $k_1!k_2!\ldots k_s!$. On the other hand, we have
\begin{eqnarray*}
(1_{G_{\alpha}})^G&=&(1_{S_{M(\alpha)}})^{S_m}\\
                  &=&\sum_{M(\alpha)\unlhd \pi}K_{\pi, M(\alpha)}\chi^{\pi}.
\end{eqnarray*}
Now, using Theorem A, we obtain,\\

{\bf Theorem B}.
$$
Q_d=\frac{1}{m!}\sum_{\alpha\in \Gamma^+_{m,d}}\sum_{M(\alpha)\unlhd \pi}M(\alpha)!K_{\pi, M(\alpha)}\chi^{\pi}.
$$
\\

As a result, we can compute the dimension of $H_d(S_m, \chi^{\pi})$, in a new fashion. We have,
\begin{eqnarray*}
\dim H_d(S_m, \chi^{\pi})&=&\chi^{\pi}(1)[\chi^{\pi}, Q_d]_{S_m}\\
                         &=&\frac{\chi^{\pi}(1)}{m!}\sum_{\alpha\in \Gamma^+_{m,d},M(\alpha)\unlhd \pi}M(\alpha)!K_{\pi, M(\alpha)}.
\end{eqnarray*}
Note that, this generalizes the similar formulae in the final part of the second section of \cite{relative}. Now, as a final result,
we have also, a necessary and sufficient condition for $H_d(S_m, \chi^{\pi})$ to be non-zero.\\

{\bf Theorem C}.\\
{\em We have $H_d(S_m, \chi^{\pi})\neq 0$, if and only if there exists $\alpha\in \Gamma^+_{m,d}$, such that $M(\alpha)\unlhd \pi$.}\\

\end{document}